\documentclass{svmult}

\renewcommand{\email}[1]{\emailname: #1} 

\usepackage{mathptmx}       
\usepackage{helvet}         
\usepackage{courier}        

\usepackage{makeidx}         
\usepackage{graphicx}        
\usepackage[bottom]{footmisc}

\usepackage{latexsym}
\usepackage{amsmath}
\usepackage{amsfonts}
\usepackage{amssymb}
\usepackage{bm}

\usepackage{url}
\usepackage{algorithm}
\usepackage{algorithmic}
\usepackage[misc,geometry]{ifsym}

\spdefaulttheorem{assumption}{Assumption}{\upshape \bfseries}{\itshape}
\spdefaulttheorem{algo}{Algorithm}{\upshape \bfseries}{\itshape}

\newcommand{\bbS}{{\mathbb{S}}}

\newcommand{\R}{{\mathbb{R}}} 

\DeclareSymbolFont{bbold}{U}{bbold}{m}{n}
\DeclareSymbolFontAlphabet{\mathbbold}{bbold}

\usepackage{amsmath,amsfonts,amssymb}
\usepackage{graphicx}

\begin{document}
\title*{A Note on Some Approximation Kernels on the Sphere}
\author{Peter Grabner}
\institute{Peter Grabner (\Letter)\at
Institut f\"ur Analysis und Zahlentheorie,
Technische Universit\"at Graz,
Kopernikusgasse 24, 8010 Graz,
Austria, \email{peter.grabner@tugraz.at}}

\maketitle

\paragraph{ Dedicated to Ian H. Sloan on the occasion of his
  80\textsuperscript{th} birthday.}

\index{Grabner, Peter} \abstract{We produce precise estimates for the
  Kogbetliantz kernel for the approximation of functions on the
  sphere. Furthermore, we propose and study a new approximation
  kernel, which has slightly better properties.} 

\date{\today}

\section{Introduction}
For $d\geq1$, let
$\bbS^d=\{\mathbf{z}\in\R^{d+1}:\langle\mathbf{z},\mathbf{z}\rangle=1\}$ denote
the $d$-dimensional unit sphere embedded in the Euclidean space $\R^{d+1}$ and
$\langle\cdot,\cdot\rangle$ be the usual inner product. We use $\mathrm{d}\sigma_d$
for the surface element and set $\omega_d=\int_{\bbS^d}\mathrm{d}\sigma_d$.

In \cite{Kogbetliantz1924:recherches_sommabilite} E.~Kogbetliantz
studied Ces\`aro means of the ultraspherical Dirichlet kernel. Let
$C_n^\lambda$ denote the $n$-th Gegenbauer polynomial of index
$\lambda$. Then for $\lambda=\frac{d-1}2$
\begin{equation*}
  K_n^{\lambda,0}(\langle\mathbf{x},\mathbf{y}\rangle)=
\sum_{k=0}^n\frac{k+\lambda}\lambda 
C_k^\lambda(\langle\mathbf{x},\mathbf{y}\rangle)
\end{equation*}
is the projection kernel on the space of harmonic polynomials of
degree $\leq n$ on the sphere $\mathbb{S}^d$. The kernel could be
studied for all $\lambda>0$, but since we have the application to polynomial
approximation on the sphere in mind, we restrict ourselves to half-integer and
integer values of $\lambda$. Throughout this paper $d$ will denote the
dimension of the sphere and $\lambda=\frac{d-1}2$ will be the
corresponding Gegenbauer parameter.

Kogbetliantz \cite{Kogbetliantz1924:recherches_sommabilite} studied
how higher Ces{\`a}ro-means improve the properties of the kernel
$K_n^{\lambda,0}$: for $\alpha\geq0$ set
\begin{equation*}
  K_n^{\lambda,\alpha}(t)=\frac1{\binom{n+\alpha}n}
\sum_{k=0}^n\binom{n-k+\alpha}{n-k}\frac{k+\lambda}\lambda C_k^\lambda(t).
\end{equation*}
He proved that the kernels $(K_n^{\lambda,\alpha})_n$ have uniformly
bounded $L^1$-norm, if $\alpha>\lambda$ and that they are
non-negative, if $\alpha\geq2\lambda+1$. There is a very short and
transparent proof of the second fact due to Reimer
\cite{Reimer1996:kogbetliantz_cesaro}. In this paper, we will
restrict our interest to the kernel $K_n^{\lambda,2\lambda+1}$, which
we will denote by $K_n^\lambda$ for short.

The purpose of this note is to improve Kogbetliantz' upper bounds for
the kernel $K_n^\lambda$. Especially, the estimates for
$K_n^\lambda(t)$ given in
\cite{Kogbetliantz1924:recherches_sommabilite} exhibit rather bad
behaviour at $t=-1$. This is partly a consequence of the actual
properties of the kernel at that point, but to some extent the
estimate used loses more than necessary. Furthermore, the estimates
given in \cite{Kogbetliantz1924:recherches_sommabilite} contain
unspecified constants. We have used some effort to provide good
explicit constants.

In the end of this paper we
will propose a slight modification of the kernel function, which is
better behaved at $t=-1$ and still shares all desirable properties of
$K_n^\lambda$.

\section{Estimating the kernel function}\label{kern}
In the following we will use the notation
$$A_n^\alpha=\binom{n+\alpha}n.$$
Notice that
\begin{equation}
\label{eq:binom}
\sum_{n=0}^\infty A_n^\alpha z^n=\frac1{(1-z)^{\alpha+1}}.
\end{equation}
Let $C_n^\lambda$ denote the $n$-th Gegenbauer polynomial with index
$\lambda$. The Gegenbauer polynomials satisfy two basic generating
function relations
(cf.~\cite{Andrews_Askey_Roy1999:special_functions,
Kogbetliantz1924:recherches_sommabilite})
\begin{align}
\label{eq:gegen1}
\sum_{n=0}^\infty C_n^\lambda(\cos\vartheta)z^n
&=\frac1{(1-2z\cos\vartheta+z^2)^\lambda}\\
\sum_{n=0}^\infty \frac{n+\lambda}\lambda C_n^\lambda(\cos\vartheta)z^n
&=\frac{1-z^2}{(1-2z\cos\vartheta+z^2)^{\lambda+1}}.\label{eq:gegen2}
\end{align}

Several different kernel functions for approximation of functions on
the sphere and their saturation behaviour have been studied in
\cite{Berens_Butzer_Pawelke1969:limitierungsverfahren_kugelfunktionen}. We
will investigate the kernel
$$K_n^\lambda(\cos\vartheta)=\frac1{A_n^{2\lambda+1}}\sum_{k=0}^n
A_{n-k}^{2\lambda+1}\frac{k+\lambda}\lambda\,C_k^\lambda(\cos\vartheta),$$
which has been shown to be positive by
E.~Kogbetliantz~\cite{Kogbetliantz1924:recherches_sommabilite} for
$\lambda>0$. 

By the generating functions \eqref{eq:binom} and \eqref{eq:gegen2} it follows
\begin{equation}
\label{eq:Kn}
\sum_{n=0}^\infty A_n^{2\lambda+1}K_n^\lambda(\cos\vartheta)z^n=
\frac{1+z}{(1-2z\cos\vartheta+z^2)^{\lambda+1}(1-z)^{2\lambda+1}}.
\end{equation}
Thus we can derive integral representations for $K_n^\lambda$ using
Cauchy's integral formula. As pointed out in the introduction, we will
restrict the values of $\lambda$ to integers or half-integers. The
main advantage of this is the fact that the exponent of $(1-z)$ in
\eqref{eq:Kn} is then an integer.

For $\lambda=k\in\mathbb{N}_0$ we split the generating function
\eqref{eq:Kn} into two factors
\begin{equation*}
  \frac{1+z}{(1-2z\cos\vartheta+z^2)(1-z)}\times
\frac1{(1-2z\cos\vartheta+z^2)^k(1-z)^{2k}}.
\end{equation*}
The first factor is essentially the generating function of the Fej\'er kernel,
namely
\begin{equation}\label{eq:fejer}
  \frac1{2\pi i}\oint\limits_{|z|=\frac12}\frac{1+z}{(1-2z\cos\vartheta+z^2)(1-z)}
\frac{\mathrm{d} z}{z^{n+1}}=
\left(\frac{\sin(n+1)\frac\vartheta2}{\sin\frac\vartheta2}\right)^2
\leq\frac1{(\sin\frac\vartheta2)^2}.
\end{equation}
Notice that this is just the kernel $(n+1)K_n^0$.

We compute the coefficients of the second factor
using Cauchy's formula
\begin{equation}
  \label{eq:cauchy}
  Q_n^k(\cos(\vartheta))=\frac1{2\pi i}\oint\limits_{|z|=\frac12}
\frac1{(1-2z\cos\vartheta+z^2)^k(1-z)^{2k}}\frac{\mathrm{d} z}{z^{n+1}}.
\end{equation}

In order to produce an estimate for $Q_n^k$, we first compute
$Q_n^1$. This is done by residue calculus and yields
\begin{equation}\label{eq:Q1}
  Q_n^1(\cos(\vartheta))=\frac1{4\sin^2(\frac\vartheta2)}
\left(n+2-\frac{\sin((n+2)\vartheta)}{\sin(\vartheta)}\right).
\end{equation}
This function is obviously non-negative and satisfies
\begin{equation}
  \label{eq:Q1-est}
  Q_n^1(\cos(\vartheta))\leq\frac{n+2}{2\sin^2(\frac\vartheta2)}.
\end{equation}

Now the functions $Q_n^k$ are formed from $Q_n^1$ by successive
convolution:
\begin{equation*}
  Q_n^{k+1}(\cos(\vartheta))=\sum_{m=0}^nQ_m^k(\cos(\vartheta))Q_{n-m}^1
(\cos(\vartheta)).
\end{equation*}
Inserting the estimate \eqref{eq:Q1-est} and an easy induction yields
\begin{equation}
  \label{eq:Q-est}
  Q_n^k(\cos(\vartheta))\leq\frac1{2^k\sin^{2k}(\frac\vartheta2)}
\sum_{r=0}^k\binom{k}r\binom{n+r+k-1}{n}.
\end{equation}
\begin{remark}\label{rem:best}
  Asymptotically, this estimate is off by a factor of $2^\lambda$, but as
  opposed to Kogbetliantz' estimate it does not contain a negative
  power of $\sin(\vartheta)$, which would blow up at $\vartheta=\pi$. The size
  of the constant is lost in the transition from \eqref{eq:Q1} to
  \eqref{eq:Q1-est}, where the trigonometric term (actually a
  Chebyshev polynomial of the second kind) is estimated by its
  maximum. On the one hand this avoids a power of $\sin(\vartheta)$ in the
  denominator, on the other hand it spoils the constant.
\end{remark}

Putting \eqref{eq:fejer} and \eqref{eq:Q-est} together yields
\begin{equation}
  \label{eq:Kn-int}
  A_n^{2k+1}K_n^k(\cos(\vartheta))\leq
\frac1{2^k(\sin\frac\vartheta2)^{2k+2}}\sum_{\ell=0}^k\binom{k}\ell
\binom{n+k+\ell}n,
\end{equation}
where we have used the identity
\begin{equation*}
  \sum_{i=0}^n\binom{i+m}i=\binom{n+m+1}n.
\end{equation*}

\begin{remark}
  Since the generating function of $A_n^{2k+1}K_n^k(\cos(\vartheta))$ is a
  rational function in this case, an application of residue calculus
  would have of course been an option. The calculation of the residues
  at $e^{\pm i\vartheta}$ produces a denominator containing
  $\sin(\vartheta)^{2k-1}$. Computation of the numerators for small values
  of $k$ show that this denominator actually cancels, but we did not
  succeed in proving this in general. Furthermore, keeping track of
  the estimates through this cancellation seems to be difficult. This
  denominator could also be eliminated by restricting $\frac
  Cn\leq\vartheta\leq\pi-\frac Cn$, but this usually spoils any gain in
  the constants obtained before. This was actually the technique used
  in \cite{Kogbetliantz1924:recherches_sommabilite}.
\end{remark}

For $\lambda=\frac12+k$ we split the generating function
\eqref{eq:Kn} into the factors
\begin{equation}\label{eq:factor}
  \frac{1}{\sqrt{1-2z\cos\vartheta+z^2}(1-z)}\times
\frac{1+z}{(1-2z\cos\vartheta+z^2)^{k+1}(1-z)^{2k+1}}
\end{equation}
with $k\in\mathbb{N}_0$. The second factor is exactly the generating
function related to the case of integer parameter $\lambda$ studied above.

For the coefficients of the first factor in \eqref{eq:factor} we use
Cauchy's formula again
\begin{equation*}
  R_n(\cos\vartheta)=\frac1{2\pi i}\oint\limits_{|z|=\frac12}
\frac{1}{\sqrt{1-2z\cos\vartheta+z^2}(1-z)}\frac{\mathrm{d} z}{z^{n+1}}.
\end{equation*}
\begin{figure}[h]
  \centering
  \includegraphics[width=7cm]{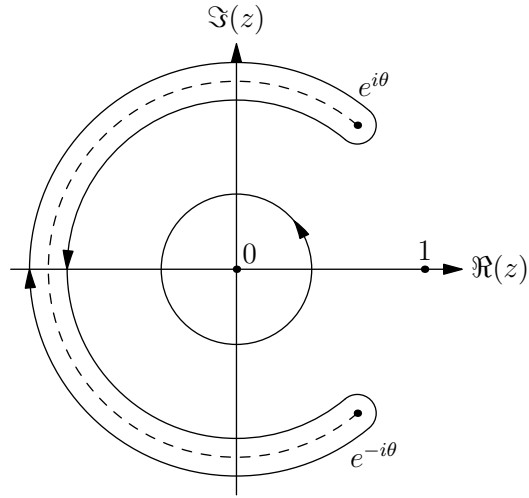}
  \caption{The contour of integration used for deriving $R_n(\cos\vartheta)$.}
  \label{fig:R_n}
\end{figure}
We deform the contour of integration to encircle the branch cut of the
square root, which is chosen to be the arc of the circle of radius one
connecting the points $e^{\pm i\vartheta}$ passing through $-1$. This
deformation of the contour passes through $\infty$ and the simple pole
at $z=1$, where we collect a residue. This gives
\begin{equation*}
  R_n(\cos\vartheta)=\frac1{2\sin\frac\vartheta2}-
\frac1{2\sqrt2 \pi}\int\limits_\vartheta^{2\pi-\vartheta}\frac{\cos((n+1)\vartheta)}
{\sqrt{\cos\vartheta-\cos t}\sin\frac t2}\,\mathrm{d} t.
\end{equation*}
We estimate this by
\begin{equation}\label{eq:R_n-est}
  R_n(\cos\vartheta)\leq\frac1{2\sin\frac\vartheta2}+\frac1{2\sqrt2\pi}
\int\limits_\vartheta^{2\pi-\vartheta}\frac{1}
{\sqrt{\cos\vartheta-\cos t}\sin\frac t2}\,\mathrm{d} t
=\frac1{\sin\frac\vartheta2}.
\end{equation}
This estimate is the best possible independent of $n$, because $R_{2n}(-1)=1$.

Putting the estimates \eqref{eq:Kn-int} and \eqref{eq:R_n-est}
together we obtain
\begin{equation}
  \label{eq:Kn-half}
  A_n^{2k+2}K_n^{k+\frac12}(\cos(\vartheta))
\leq
\frac1{2^{k}(\sin(\frac\vartheta2))^{2k+3}}\sum_{\ell=0}^{k}\binom{k}\ell
\binom{n+k+\ell+1}n.
\end{equation}

Summing up, we have proved the following.
\begin{theorem}\label{thm:kogb} Let $\lambda=\frac{d-1}2$ be a
  positive integer or half-integer. 
  Then the kernel $K_n^\lambda$ satisfies the following estimates
  \begin{equation}
    \label{eq:K_n-est}
    K_n^\lambda(\cos\vartheta)\leq
    \begin{cases}
\displaystyle{\frac1{2^{\lfloor\lambda\rfloor}(\sin(\frac\vartheta2))^{2\lambda+2}}
\sum_{\ell=0}^{\lfloor\lambda\rfloor}\binom{\lfloor\lambda\rfloor}\ell
\frac{(2\lambda+1)_{\ell+1}}
{(n+2\lambda+1)_{\ell+1}}}&\text{for}\quad0<\vartheta\leq\pi\\[6mm]
\frac{(n+4\lambda+1)_n}{(n+2\lambda)_n}&\text{for}\quad0\leq\vartheta\leq\pi,
    \end{cases}
  \end{equation}
where $(a)_n=a(a-1)\cdots(a-n+1)$ denotes the \textbf{falling}
factorial (Pochhammer symbol).
\end{theorem}
\begin{remark}
  The estimate \eqref{eq:K_n-est} is best possible with respect to the
  behaviour in $n$ for a fixed $\vartheta\in(0,\pi)$, as well as for the
  power of $\sin\frac\vartheta2$. The constant in front of the main
  asymptotic term could still be improved, especially its dependence
  on the dimension. The second estimate is the trivial estimate by
  $K_n^\lambda(1)$.
\end{remark}
\section{A new kernel}
\begin{figure}[t]
  \centering
  \includegraphics[width=\textwidth]{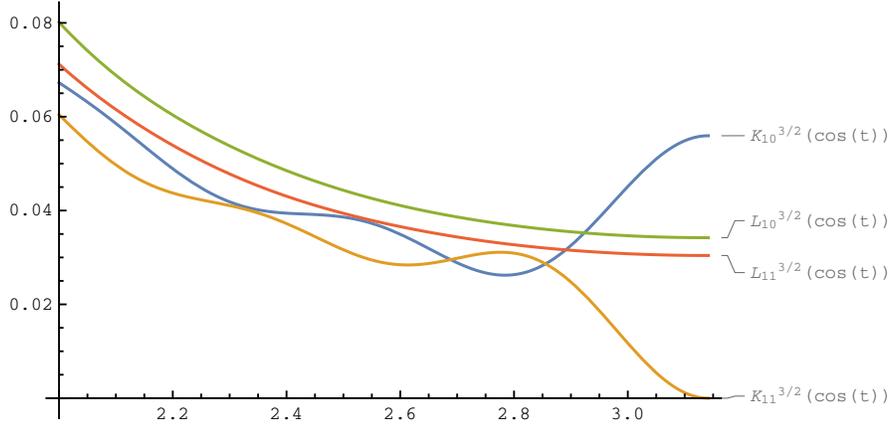}
  \caption{Comparison between the kernels $K_{10}^{\frac32}$,
    $K_{11}^{\frac32}$, $L_{10}^{\frac32}$, and $L_{11}^{\frac32}$. The
    kernels $K$ show oscillations and a parity phenomenon at $\vartheta=\pi$.}
  \label{fig:comparison}
\end{figure}

The kernel $K_n^\lambda(\cos\vartheta)$ exhibits a parity phenomenon at
$\vartheta=\pi$, which occurs in the first asymptotic order term (see
Figure~\ref{fig:comparison} for illustration). This
comes from the fact that the two singularities at $e^{\pm i\vartheta}$
collapse to one singularity of twice the original order for this value
of $\vartheta$. In order to avoid this, we propose to study the kernel
given by the generating function
\begin{equation}
  \label{eq:L-gen}
  \frac{(1+z)^{2\lambda+2}}{(1-2z\cos\vartheta+z^2)^{\lambda+1}(1-z)^{2\lambda+1}}=
\frac{1-z^2}{(1-2z\cos\vartheta+z^2)^{\lambda+1}}\times
\frac{(1+z)^{2\lambda+1}}{(1-z)^{2\lambda+2}}.
\end{equation}
Let $B_n^\lambda$ be given by
\begin{equation}\label{eq:Bn-generating}
  \sum_{n=0}^\infty B_n^\lambda z^n=\frac{(1+z)^{2\lambda+1}}{(1-z)^{2\lambda+2}},
\end{equation}
then the kernel is given by
\begin{align}
  \label{eq:L-def}
  L_n^\lambda(\cos\vartheta)&=\frac1{B_n^\lambda}\sum_{k=0}^nB_{n-k}^\lambda 
\frac{k+\lambda}\lambda C_k^\lambda(\cos\vartheta)\\
&=\frac1{B_n^\lambda}\sum_{\ell=0}^{2\lambda+1}\binom{2\lambda+1}\ell 
A_{n-\ell}^{2\lambda+1}K_{n-\ell}^\lambda(\cos\vartheta).
\end{align}
The coefficients $B_n^\lambda$ satisfy
\begin{align*}
  B_n^\lambda&=\sum_{\ell=0}^{2\lambda+1}\binom{2\lambda+1}\ell 
\binom{n-\ell+2\lambda+1}{n-\ell}\\&
=\sum_{\ell=0}^{2\lambda+1}(-1)^\ell\binom{2\lambda+1}\ell2^{2\lambda+1-\ell}
\binom{n-\ell+2\lambda+1}{n}\sim \frac{2^{2\lambda+1}n^{2\lambda+1}}
{(2\lambda+1)!}.
\end{align*}
The expression in the second line, which allows to read of the
asymptotic behaviour immediately, is obtained by expanding the
numerator in \eqref{eq:Bn-generating} into powers of $1-z$.

For $\lambda\in\mathbb{N}_0$ we write the generating function of
$B_n^\lambda L_n^\lambda(\cos\vartheta)$ as
\begin{equation}
  \label{eq:factor-L}
  \left(\frac{(1+z)^2}{(1-2z\cos\vartheta+z^2)(1-z)^2}\right)^\lambda\times
\frac{(1+z)^2}{(1-2z\cos\vartheta+z^2)(1-z)}.
\end{equation}
The coefficients of the first factor are denoted by
$S_n^\lambda(\cos\vartheta)$. They are obtained by successive convolution of
\begin{align*}
S_n^1(\cos\vartheta)&= 
\frac1{2\pi i}\oint\limits_{|z|=\frac12}\frac{(1+z)^2}
{(1-2z\cos\vartheta+z^2)(1-z)^2}\frac{\mathrm{d} z}{z^{n+1}}\\
&=\frac{n+1}{\sin^2\frac\vartheta2}\left(1-
 \frac{\cos (\frac\vartheta2 )\sin (n+1) \vartheta} 
{2(n+1)\sin\frac{\vartheta }{2}}\right).
\end{align*}
In order to estimate $S_n^1(\cos\vartheta)$, we estimate the
$\mathrm{sinc}$-function by its minimum
\begin{equation*}
\mathrm{sinc}(t)=\frac{\sin(t)}t\geq -C'=-0.217233628211221657408279325562\ldots.
\end{equation*}
The value was obtained with the help of \texttt{Mathematica}. This gives
\begin{align*}
  &1-\cos\left(\frac\vartheta2\right)\frac{\sin((n+1)\vartheta)}{2(n+1)\sin(\frac\vartheta2)}
=1-\mathrm{sinc}((n+1)\vartheta)\frac{\cos\frac\vartheta2}{\mathrm{sinc}\,\frac\vartheta2}\\
&\leq1+C'=:C=1.217233628211221657408279325562\ldots,
\end{align*}
where we have used that $\cos(\frac\vartheta2)\leq\mathrm{sinc}(\frac\vartheta2)$ for
$0\leq\vartheta\leq\pi$.
From this we get the estimate
\begin{equation*}
  S_n^1(\cos\vartheta)\leq C\frac{n+1}{\sin^2\frac\vartheta2}
\end{equation*}
and consequently
\begin{equation}
  \label{eq:S-est}
  S_n^\lambda(\cos\vartheta)\leq \frac{C^\lambda}{\sin^{2\lambda}\frac\vartheta2} 
\binom{n+2\lambda-1}n
\end{equation}
by successive convolution as before.
\begin{remark}
  This expression is bit simpler than the corresponding estimate for
  $Q_n^\lambda$, because the iterated convolution of the terms $n+1$
  is a binomial coefficient, whereas the iterated convolution of terms
  $n+2$ can only be expressed as a linear combination of binomial
  coefficients. The growth order is the same.
\end{remark}

In a similar way we estimate the coefficient of the second factor in
\eqref{eq:factor-L} 
\begin{align*}
&\frac1{2\pi i}\oint\limits_{|z|=\frac12}\frac{(1+z)^2}{(1-2z\cos\vartheta+z^2)(1-z)}
\frac{\mathrm{d} z}{z^{n+1}}\\
&=\frac{1}{2\sin^2\frac\vartheta2}\left(2-\cos (n \vartheta )-\cos((n+1)\vartheta)\right)
\leq \frac2{\sin^2\frac\vartheta2}.
\end{align*}
As before, this is the kernel function for $\lambda=0$.

Putting this estimate together with \eqref{eq:S-est} we obtain
\begin{equation}
  \label{eq:L-est}
  B_n^\lambda L_n^\lambda(\cos\vartheta)\leq
\frac{2C^\lambda}{\sin^{2\lambda+2}\frac\vartheta2}\binom{n+2\lambda}n
\end{equation}
for $\lambda\in\mathbb{N}_0$.

For $\lambda=k+\frac12$ ($k\in\mathbb{N}_0$) we factor the generating
function as
\begin{equation}
  \label{eq:factor-half}
  \frac{(1+z)}{\sqrt{1-2z\cos\vartheta+z^2}(1-z)}\times
\frac{(1+z)^{2k+2}}{(1-2z\cos\vartheta+z^2)^{k+1}(1-z)^{2k+1}}.
\end{equation}
We still have to estimate the coefficient of the first factor, which
is given by the integral
\begin{equation*}
  T_n(\cos\vartheta)=\frac1{2\pi i}\oint\limits_{|z|=\frac12}\frac{(1+z)}
{\sqrt{1-2z\cos\vartheta+z^2}(1-z)}\frac{\mathrm{d} z}{z^{n+1}}.
\end{equation*}
We transform this integral in the same way as we did before using the
contour in Figure~\ref{fig:R_n} which yields
\begin{equation}
  \label{eq:T_n}
  T_n(\cos\vartheta)=\frac1{\sin\frac\vartheta2}-
\frac1{\pi\sqrt2}\int\limits_\vartheta^{2\pi-\vartheta}
\frac{\cos(\frac t2)\cos((n+\frac12)t)}
{\sqrt{\cos\vartheta-\cos t}\sin\frac t2}\,\mathrm{d} t.
\end{equation}
The modulus of the integral can be estimated by
\begin{equation*}
  \frac{\sqrt2}{\pi}\int\limits_{\vartheta}^{\pi}
\frac{\cos(\frac t2)}
{\sqrt{\cos\vartheta-\cos t}\sin\frac t2}\,\mathrm{d}
t=\frac{\pi-\vartheta}{\pi\sin\frac\vartheta2}\leq \frac1{\sin\frac\vartheta2}.
\end{equation*}
This gives the bound
\begin{equation}
  \label{eq:T_n-est}
  T_n(\cos\vartheta)\leq \frac2{\sin\frac\vartheta2}.
\end{equation}
Putting this estimate together with \eqref{eq:L-est} we obtain
\begin{equation}
  \label{eq:L-est-half}
  B_n^\lambda L_n^\lambda(\cos\vartheta)\leq
  \frac{4C^k}{\sin^{2k+3}\frac\vartheta2}
\binom{n+2k+1}{n}
\end{equation}
for $\lambda=k+\frac12$.

Summing up, we have proved the following. As before, the second
estimate is just the trivial estimate by $L_n^\lambda(1)$.
\begin{theorem}\label{thm:kogb} Let $\lambda=\frac{d-1}2$ be a
  positive integer or half-integer. 
  Then the kernel $L_n^\lambda$ satisfies the following estimates
  \begin{equation}
    \label{eq:L_n-est}
    L_n^\lambda(\cos\vartheta)\leq
    \begin{cases}
      D_\lambda\frac{C^{\lfloor\lambda\rfloor}}
{B_n^\lambda \sin^{2\lambda+2}\frac\vartheta2}\binom{n+2\lambda}n
&\text{for}\quad0<\vartheta\leq\pi\\[6mm]
\frac1{B_n^\lambda}
\sum_{\ell=0}^{2\lambda+2}\binom{2\lambda+2}\ell2^{2\lambda+2-\ell}(-1)^\ell
\binom{n+4\lambda+2-\ell}n&\text{for}\quad0\leq\vartheta\leq\pi,
    \end{cases}
  \end{equation}
where $D_\lambda=2$ for $\lambda\in\mathbb{N}$ and
$D_\lambda=4$, if $\lambda\in\frac12+\mathbb{N}_0$.
\end{theorem}
\begin{figure}[h]\label{fig:order}
  \centering
\includegraphics[width=\hsize]{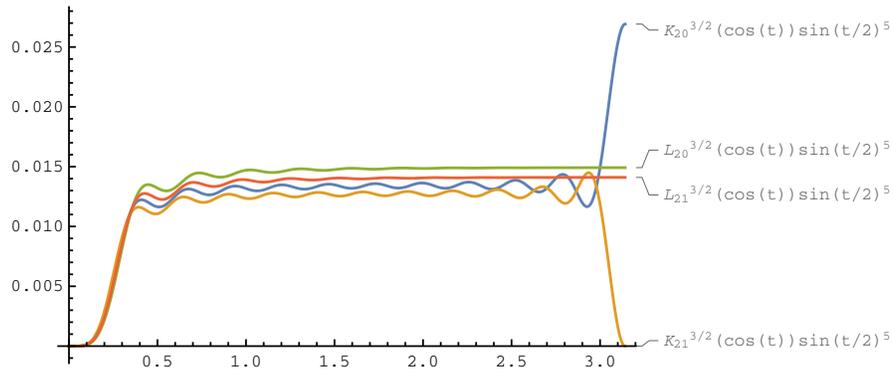}
\caption{Plots of the functions
  $K_{20}^{\frac32}(\cos\vartheta)(\sin\frac\vartheta2)^5$,
  $K_{21}^{\frac32}(\cos\vartheta)(\sin\frac\vartheta2)^5$,
  $L_{20}^{\frac32}(\cos\vartheta)(\sin\frac\vartheta2)^5$, and
  $L_{21}^{\frac32}(\cos\vartheta)(\sin\frac\vartheta2)^5$. Again the parity
  phenomenon for the kernel $K$ is prominently visible.}
\end{figure}
\begin{remark}
  Notice that the orders of magnitude in terms of $n$ and the powers
  of $\sin\frac\vartheta2$ are the same for $L_n^\lambda$ as for the
  kernel $K_n^\lambda$. This fact is illustrated by
  Figure~\ref{fig:order}. The coefficient of the asymptotic leading
  term of the estimate decays like $(2\lambda+1)(C/4)^\lambda$ for
  $L_n^\lambda$, whereas this coefficient decays like
  $(2\lambda+1)(1/2)^\lambda$ for $K_n^\lambda$.
\end{remark}

\begin{acknowledgement}
  The author is supported by the Austrian Science Fund FWF projects
  F5503 (part of the Special Research Program (SFB) ``Quasi-Monte
  Carlo Methods: Theory and Applications'') and W1230 (Doctoral
  Program ``Discrete Mathematics''). The author is grateful to two
  anonymous referees for their many helpful comments.
\end{acknowledgement}

\end{document}